# MAXIMUM FISHER INFORMATION IN MIXED STATE QUANTUM SYSTEMS

By Alessandra Luati

*University of Bologna*

We deal with the maximization of classical Fisher information in a quantum system depending on an unknown parameter. This problem has been raised by physicists, who defined [Helstrom (1967) *Phys. Lett. A* **25** 101–102] a quantum counterpart of classical Fisher information, which has been found to constitute an upper bound for classical information itself [Braunstein and Caves (1994) *Phys. Rev. Lett.* **72** 3439–3443]. It has then become of relevant interest among statisticians, who investigated the relations between classical and quantum information and derived a condition for equality in the particular case of two-dimensional pure state systems [Barndorff-Nielsen and Gill (2000) *J. Phys. A* **33** 4481–4490].

In this paper we show that this condition holds even in the more general setting of two-dimensional mixed state systems. We also derive the expression of the maximum Fisher information achievable and its relation with that attainable in pure states.

**1. Introduction.** Quantum statistics is ordinary statistical inference applied to quantum systems. The methodology is based on the mathematical specification of the state of the quantum system, to be denoted by $\rho = \rho(\theta)$ as it is supposed to depend on an unknown parameter $\theta$, and of the measurement $M$ to be carried out on that system. In finite-dimensional quantum systems, both $\rho(\theta)$ and $M$ are represented by Hermitian matrices. With such mathematical specifications, we will be able to compute the probability distribution of a random variable $X$ when a measurement $M$ is carried out on the system in state $\rho$, that is, $P_X(\cdot, \theta) = \text{tr}\{\rho(\theta)M(\cdot)\}$. As $\rho(\theta)$ depends on the parameter $\theta$, $P_X(\cdot, \theta)$ will depend on $\theta$ too, thereby setting the statistical problem of how best to estimate the value of $\theta$. Since one often has a choice of what measurement to take, the design problem of how to best measure $\rho(\theta)$ arises. In other words, the question is which measurement provides more statistical information about the unknown parameter $\theta$.

By "statistical information" we mean expected Fisher information $i(\theta, M)$, a measure of how precise an unbiased estimator $t(x)$ of $\theta$ based on the outcome of an arbitrary measurement $M$ is, as follows by the Cramér–Rao bound on the variance $V$ of $t(x)$,

$$V\{t(x)\} \geq i(\theta, M)^{-1}.$$

Based on an operator called "symmetric logarithmic derivative," which is a noncommutative logarithmic derivative for matrices, Helstrom (1967) derived a quantity, $I(\theta)$, that he called "quantum expected information" because of the relation

$$V\{t(x)\} \geq I(\theta)^{-1},$$

known as the quantum Cramér–Rao bound. The proof [see also Holevo (1982)] essentially follows the derivation of the classical Cramér–Rao bound.

Braunstein and Caves (1994) have emphasized the relation between classical and quantum information, deriving the *information inequality*

$$i(\theta, M) \leq I(\theta) \tag{1}$$

from which one obtains Helstrom's bound as a corollary. The proof of the information inequality for the one-dimensional parameter case is based on the Cauchy–Schwarz inequality with the Hilbert–Schmidt inner product and gives as a by-product two equality conditions on the measurement $M$. So, in some sense, for a particular value of the parameter, one can say what is the best measurement in terms of Fisher information. However, the Braunstein and Caves equality conditions are not very transparent. In particular, it is not clear which characteristics a measurement has to have in order to allow attainability in the information inequality. The same problem is considered by Barndorff-Nielsen and Gill (2000) dealing with quantum information. The contents of the paper are rich and the solution is elegant; specializing to the most simple possible case, the so-called pure state, two-dimensional, the authors show that a necessary and sufficient condition for attainability is that each measurement is proportional to a rank-one projection matrix.

In this paper, we consider the more general case of mixed states, and we show that the necessary and sufficient condition for pure states holds as well. The key argument in our proof is the derivation of the symmetric logarithmic derivative of a mixed state. We also derive the expression of the maximum achievable Fisher information.

FISHER INFORMATION IN QUANTUM SYSTEMS 3The paper is organized as follows. In Section 2 we introduce the basic concepts of quantum statistical inference and review the attainability condition for two-dimensional pure states. In Section 3 we generalize this result to two-dimensional mixed states. A discussion is in Section 4 and Section 5 contains the proofs of the main results.

**2. Quantum statistics and the information inequality.** In quantum statistics [Holevo (1982) and Helstrom (1976)], the probability distribution of a random variable $X : (\Omega, \mathcal{F}, P) \to (\mathbb{G}, \mathcal{G}, P_X)$ is given by the *trace rule for probability*

$$P_X(G; \theta) = \text{tr}\{\rho(\theta) M(G)\} \qquad \forall G \in \mathcal{G},$$

where $\rho(\theta)$ is a *density matrix*, that is, a nonnegative, self-adjoint and trace-one linear operator acting on a $n$-dimensional complex Hilbert space $\mathcal{H}_n$ and depending on an unknown parameter $\theta \in \Theta \subseteq \mathbb{R}^k$, while $M$ is an *operator-valued probability measure*, that is, a set of nonnegative and self-adjoint linear operators defined on the measure space $(\mathbb{G}, \mathcal{G})$ and taking values in $\mathcal{H}_n$, such that $M(\mathbb{G}) = \mathsf{I}$, the identity operator, $M(\varnothing) = \mathsf{O}$, the null operator, and

$$M\left(\bigcup_{h=1}^{\infty} G_h\right) = \sum_{h=1}^{\infty} M(G_h), \qquad \text{if } G = \bigcup_{h=1}^{\infty} G_h, G_h \cap G_k = \varnothing,$$

$\forall h, k = 1, \ldots, \infty, h \neq k$. If the measurement $M$ is absolutely continuous with respect to a $\sigma$-finite measure $\mu$ on $(\mathbb{G}, \mathcal{G})$ such that $M(G) = \int_G m(x) \mu(dx)$ $\forall G$, where $m(x)$ is nonnegative and Hermitian, then $\{m(x)\}_{x \in \mathbb{G}}$ is a *generalized measurement* and the probability density of $X$ is

$$p(x; \theta) = \text{tr}\{\rho(\theta) m(x)\}.$$

If $n < \infty$, then $\mathcal{H}_n$ can be identified with the $n$-dimensional Euclidean complex space $\mathbb{C}^n$ and it is equivalent to talk about self-adjoint operators or Hermitian matrices.

Once a *parametric quantum model* $\{\rho(\theta), M; \theta \in \Theta \subseteq \mathbb{R}^k\}$ has been chosen to describe the set of probabilistic outcomes of a random experiment, the *expected Fisher information* can be obtained as

$$i(\theta, M) = E\{l_{/\theta}^2\} = \int_{\mathbb{G}_+} p(x; \theta)^{-1} \text{tr}^2\{\rho_{/\theta}(\theta) m(x)\} \mu(dx),$$

where $l_{/\theta}$ is the score function of $\theta$, $\mathbb{G}_+ = \{x \in \mathbb{G} : p(x; \theta) > 0\}$, and $\rho_{/\theta}$ is the matrix whose $ij$th generic element is $[\rho_{/\theta}]_{ij} = \frac{\partial}{\partial \theta}[\rho(\theta)]_{ij}$.

Braunstein and Caves (1994) showed that at a fixed value of the parameter $\theta$, an upper bound for Fisher information is the *quantum information*

(2) $$I(\theta) = E\{\rho_{//\theta}^2\} = \text{tr}\{\rho(\theta) \rho_{//\theta}^2\}$$

introduced by Helstrom (1967) based on the *symmetric logarithmic derivative* (SLD) or *symmetric quantum score* $\rho_{//\theta}$, implicitly defined by the relation

(3) $$\rho_{/\theta} = \tfrac{1}{2}[\rho(\theta)\rho_{//\theta} + \rho_{//\theta}\rho(\theta)].$$

Barndorff-Nielsen and Gill (2000) showed that equality holds in (1) if

(4) $$k(x,\theta)^{1/2} m(x)^{1/2} \rho(\theta)^{1/2} = m(x)^{1/2} \rho_{//\theta} \rho(\theta)^{1/2}$$

for $\mu$ almost all $x$ in $\mathbb{G}_+$, where $k(x,\theta) = p(x;\theta)^{-1} \operatorname{tr}\{\rho(\theta)\rho_{//\theta} m(x) \rho_{//\theta}\}$ is real, since $k(x,\theta) = p(x;\theta)^{-1} \operatorname{tr}\{C^H C\}$ with $C = \rho^{1/2}(\theta)\rho_{//\theta} m^{1/2}(x)$ and $k(x,\theta) = 0$ over $\mathbb{G}_0 = \{x \in \mathbb{G} : p(x;\theta) = 0\}$.

We will refer to (4) as the *equality condition* between classical and quantum information and we call a measurement that satisfies it an *attaining measurement*. With some abuse of notation the single matrices that constitute the family of the attaining measurements will be called attaining measurements as well.

A characterization of an attaining measurement is given by Barndorff-Nielsen and Gill (2000) for the special case of one-parameter pure states,

$$\rho(\theta) = |\psi(\theta)\rangle\langle\psi(\theta)|,$$

where $\psi(\theta)$ is a unit vector and $\theta \in \Theta \subseteq \mathbb{R}$; in the Dirac *bra-ket* notation $|\psi(\theta)\rangle$ denotes a column vector (*ket*) and $\langle\psi(\theta)|$ its Hermitian transpose (row, *bra*). According to quantum theory, pure states represent the best knowledge one can have about some specific properties of the system under observation. The authors show that, in two-dimensional pure states, a necessary and sufficient condition for attainability is that the attaining measurement be proportional to a rank-one projection matrix. The proof is based on the following properties of pure state density matrices; in particular, expression (i) for the symmetric logarithmic derivative of a pure state plays a crucial role:

(i) $\rho_{//\theta} = 2\rho_{/\theta}$;
(ii) $\rho(\theta)\rho_{/\theta}\rho(\theta) = \mathsf{O}$;
(iii) $\operatorname{tr}\{\rho_{//\theta}\rho(\theta)\} = 0$;
(iv) $I(\theta) = 2\operatorname{tr}\{\rho_{/\theta}^2\}$.

These properties are no more than algebraic consequences of the definitions of pure state and symmetric logarithmic derivative. An interesting way of proving them is in Fujiwara and Nagaoka (1995) where they are derived as properties of two preinner products defined on the set of all the linear and bounded operators on $\mathcal{H}$, in the context of estimation in pure state models.

In the following, we generalize Barndorff-Nielsen and Gill's condition to the more general setting of mixed states. We derive some properties analogous to (i)–(iv) as well as the symmetric logarithmic derivative of a mixed state. This latter is the key result of the paper, since, based upon it, we obtain maximum Fisher information and the necessary and sufficient condition for a measurement to attain it.

**3. Attainability conditions in mixed states.** A quantum system is said to be a *mixed state* if its density matrix is of the form

$$\rho(\theta) = w_1(\theta)\rho_1(\theta) + w_2(\theta)\rho_2(\theta) + \cdots + w_m(\theta)\rho_m(\theta),$$

where $\rho_i(\theta) = |\psi_i(\theta)\rangle\langle\psi_i(\theta)|$, $i = 1,\ldots,m$ and $|\psi_1(\theta)\rangle, |\psi_2(\theta)\rangle, \ldots, |\psi_m(\theta)\rangle$ are unit vectors; the $w_i(\theta)$'s are real weights satisfying $w_i(\theta) \geq 0 \ \forall i = 1,\ldots,m$ and $\sum_{i=1}^{m} w_i(\theta) = 1$.

Mixed states, obtained as convex combinations of pure states, indicate a situation of partial knowledge of the system. They represent probabilistic mixtures, in the sense that the system under observation is in the state $\rho_i(\theta)$ with probability $w_i(\theta) \ \forall i = 1,\ldots,m$.

In considering mixed states, the problem of characterizing attaining measurements becomes much more complicated than in pure states. However, restricting to mixtures of two-dimensional orthogonal pure states, analogous conclusions to those for pure states can be obtained. The dimensionality constraint is not restrictive, since two-dimensional systems are the most frequently encountered in quantum mechanics, both because they are simple to deal with and for their objective importance. For instance, electrons, qubits and spin-$\frac{1}{2}$ particles are just some examples of two-dimensional systems playing a crucial role in quantum mechanics. The orthogonality constraint is one among an infinity of choices of how a mixed state can be decomposed into pure states. There is an illuminating geometrical illustration for the decomposability of mixed states which is based on the Bloch or Poincaré or Riemann sphere representation of the set of states in two-dimensional complex Hilbert spaces by means of unit vectors in real three-dimensional Euclidean spaces [see Luati (2003) and references therein]. In particular, there exists a one-to-one correspondence between states in $\mathbb{C}^2$ and the unit ball in $\mathbb{R}^3$. So, if $\mathcal{H} = \mathbb{C}^2$, then the set of pure states is the surface of the unit sphere and the set of mixed states is the interior of the corresponding unit ball. Mixtures of two pure states can be represented as points in the interior of the sphere, on the straight line joining the two points on the surface. If the generating pure states are orthogonal (opposite on the sphere), then the corresponding mixed states lie on the diameters of the great circles. Therefore, the set of such states with given weights can be represented by the spheres embedded in the unit sphere with the same center, but radius smaller than one and dependent on the weights of the mixtures. As we will see, this characteristic plays a relevant role in the geometric interpretation of the results that follow.

A one-parameter two-dimensional mixed state can be represented as

$$\rho(\theta) = w(\theta)\rho_1(\theta) + (1 - w(\theta))\rho_2(\theta), \tag{5}$$

where $\rho_1(\theta) = |\psi_1(\theta)\rangle\langle\psi_1(\theta)|$ and $\rho_2(\theta) = |\psi_2(\theta)\rangle\langle\psi_2(\theta)|$ such that $\langle\psi_1(\theta)|\psi_1(\theta)\rangle = 1$ and $|\psi_2(\theta)\rangle = I_1^{-1/2}(\theta)\rho_{1//\theta}|\psi_1(\theta)\rangle$ is a unit vector with $I_1^{-1/2}(\theta)$ normalizing constant; $I_h(\theta)$ indicates quantum information provided by $\rho_h(\theta)$, and $\rho_{h/\theta}$ and $\rho_{h//\theta}$, respectively, stand for the term-by-term first derivative and for the symmetric logarithmic derivative of $\rho_h(\theta), h = 1, 2$, with respect to $\theta \in \Theta \subseteq \mathbb{R}$. We assume that the coefficient $w(\theta)$ is a smooth function of $\theta$ taking values in the real interval $(0, 1)$; we also consider $w(\theta) \neq \frac{1}{2}$ because otherwise $\rho(\theta) = \frac{1}{2}\mathsf{I}$ does not depend on any unknown parameter. This case (the center of the unit sphere) represents the maximum entropy situation, that is, complete ignorance about the quantum system under observation. The vectors $|\psi_1(\theta)\rangle$ and $|\psi_2(\theta)\rangle$ are orthonormal by (iii) and by definition of $I_1(\theta)$. It therefore follows (the arguments are sometimes omitted) that

(v) $\rho_2(\theta) = \mathsf{I} - \rho_1(\theta)$;
(vi) $\rho_{2/\theta} = -\rho_{1/\theta}$;
(vii) $\rho_{2//\theta} = -\rho_{1//\theta}$;
(viii) $\rho_h(\theta)\rho_{k/\theta}\rho_h(\theta) = \mathsf{O}$, for $h, k = 1, 2$;
(ix) $I_2(\theta) = I_1(\theta)$.

In fact, (v) is a consequence of the spectral theorem in $\mathbb{C}^2$; (vi) is obtained by differentiating term-by-term the elements of (v) with respect to $\theta$; (vii) follows by (i) for pure states and by (vi); (viii) follows by (ii) for pure states and by (vi). Finally, (ix) follows by (iv) and (vi).

As in the pure state case, to draw some conclusion on attaining measurements starting from equality condition (4) it is necessary to know the exact form of $\rho_{//\theta}$.

LEMMA 1. *The symmetric logarithmic derivative of the mixed state* (5) *is*

$$\rho_{//\theta} = \frac{w_{/\theta}}{w(\theta)}\rho_1(\theta) + (2w(\theta) - 1)\rho_{1//\theta} - \frac{w_{/\theta}}{1 - w(\theta)}\rho_2(\theta), \tag{6}$$

*where $w_{/\theta}$ stands for the first derivative (scalar) of $w(\theta)$ with respect to $\theta$.*

Based on Lemma 1, we obtain quantum or Helstrom information, that is, an upper bound for Fisher information in the mixed state (5).

LEMMA 2. *Quantum information given by the mixed state* (5) *is*

$$I(\theta) = \frac{(w_{/\theta})^2}{w(\theta)(1 - w(\theta))} + (2w(\theta) - 1)^2 I_1(\theta). \tag{7}$$

Observe that if $w$ does not depend on $\theta$, then

(8) $$I(\theta) = (2w-1)^2 I_1(\theta)$$

and, since $w \in (0,1)$, quantum information provided by a mixed state is less than quantum information provided by a pure state. This is a quantum information based way to state that pure states represent the best knowledge that one can have about a quantum system. On the other hand, if $w(\theta)$ depends on $\theta$, then no conclusions about quantum information can be drawn without knowing the function $w(\theta)$, except that $I(\theta) < I_1(\theta)$ if $\frac{(w_{/\theta})^2}{4w(\theta)^2(1-w(\theta))^2} < I_1(\theta)$. However, in both cases, we can specify the measurements such that Fisher information is maximum.

THEOREM 1. *In mixed states of the form* (5), $i(\theta, M) = I(\theta)$ *if and only if, for $\mu$-almost all $x$ in $\mathbb{G}_+$,*

$$\{m(x)\}_{x \in \mathbb{G}_+} \propto_{\mathbb{R}} \{|\gamma(x)\rangle\langle\gamma(x)|\}_{x \in \mathbb{G}_+}$$

*with*

(9) $$\langle\gamma(x)| |\psi_1(\theta)\rangle \propto_{\mathbb{R}} \langle\gamma(x)| |\psi_2(\theta)\rangle,$$

*where $\langle\gamma(x)| |\gamma(x)\rangle = 1$ and $\propto_{\mathbb{R}}$ stands for "proportional by means of a real constant."*

What is remarkable in this result is that in two-dimensional systems, pure or mixed, attaining measurements are of the same form. Furthermore, when $w$ does not depend on $\theta$, the maximum Fisher information achievable cannot be greater than that achievable in pure states. In this case, even the geometric aspects pointed out by Barndorff-Nielsen and Gill (2000), based on the sphere representation of pure states, can be interpreted in terms of mixed states as well.

In fact, these authors show that in pure states like $\rho_1(\theta)$, condition (9) geometrically implies that the attaining measurements are proportional to rank-one orthogonal projection matrices onto state vectors that correspond to points on the intersection of the unit sphere with a plane spanned by two orthogonal vectors of $\mathbb{R}^3$. These are the unit vector $\underline{u}$ such that $\rho_1(\theta) = \frac{1}{2}(\mathsf{I} + \langle \underline{u}(\theta), \underline{\sigma}\rangle)$ and the orthogonal vector $\underline{u}_{/\theta}$ such that $I_1(\theta) = \|\underline{u}_{/\theta}\|^2 > 0$, where $\underline{\sigma}^T = (\sigma_x, \sigma_y, \sigma_z)$ is the vector of Pauli matrices

$$\sigma_x = \begin{bmatrix} 0 & 1 \\ 1 & 0 \end{bmatrix}, \qquad \sigma_y = \begin{bmatrix} 0 & -i \\ i & 0 \end{bmatrix}, \qquad \sigma_z = \begin{bmatrix} 1 & 0 \\ 0 & -1 \end{bmatrix}$$

that together with the identity matrix $\mathsf{I}$ constitute an orthogonal basis for the set of Hermitian matrices acting on $\mathbb{C}^2$. Moreover, if quantum information is positive, then there exist uniformly attaining measurements, that is, such that $i(\theta; M) = I(\theta) \ \forall \theta \in \Theta$, if and only if the set of the states of the given model is a great circle in the unit sphere. Otherwise, no measurement exists such that equality holds for all the values of the parameter $\theta$.

Consider now mixed states like (5). Given $\rho_1(\theta) = \frac{1}{2}(\mathsf{I} + \langle \underline{u}(\theta), \underline{\boldsymbol{\sigma}} \rangle)$, then $\rho_2(\theta) = \frac{1}{2}(\mathsf{I} + \langle -\underline{u}(\theta), \underline{\boldsymbol{\sigma}} \rangle)$ and therefore $\rho(\theta) = \frac{1}{2}(\mathsf{I} + \langle (2w(\theta) - 1)\underline{u}(\theta), \underline{\boldsymbol{\sigma}} \rangle)$. This means that the set of mixed states with given $w(\theta)$ is a sphere of radius $|2w(\theta) - 1| < 1$ embedded in the unit ball, with the same center. Hence, when $w$ does not depend on $\theta$, up to the factor $2w - 1$ results for mixed states can be read as results for pure states. Particularly, since $\text{span}\{(2w-1)\underline{u}(\theta), \underline{u}_{/\theta}\} \equiv \text{span}\{\underline{u}(\theta), \underline{u}_{/\theta}\}$, condition (9) implies that attaining measurements are proportional to rank-one projectors onto state vectors that correspond to points on the intersection of the sphere of radius $|2w - 1|$ with the plane spanned by the vectors $\underline{u}(\theta)$ and $\underline{u}_{/\theta}$. Furthermore, $I(\theta) = (2w-1)^2 \|\underline{u}_{/\theta}\|^2 > 0$ and we can conclude that uniformly attaining measurements in mixed states are admitted if and only if the model is a great circle in the sphere of radius $|2w - 1|$.

**4. Discussion.** We characterized the measurements that maximize Fisher information in two-dimensional mixed states. This result generalizes that of Barndorff-Nielsen and Gill (2000) for pure states. We also derived the maximum Fisher information achievable in mixed states and the conditions such that it is not (or cannot be) greater than that attainable in pure states.

A particular decomposition of a two-dimensional mixed state into orthogonal pure states allowed a geometric interpretation of the results that gave special emphasis to measurements that do not depend on the unknown parameter. The interest in uniformly attaining measurements traces its origins to the same paper in which the information inequality is derived. In analogy with the metric properties of Fisher information, Braunstein and Caves proposed to use quantum information as a metric on the set of all the possible states of a given quantum system. However, this is sensible only when uniform attainability holds. It is straightforward to show that in mixed states represented through vectors parametrized by colatitude $\eta$ (known) and longitude $\phi$ in a sphere of radius $|2w - 1|$, where $\langle \psi_1(\eta, \phi)| = [\cos(\frac{\eta}{2})e^{i\phi/2} \sin(\frac{\eta}{2})e^{-i\phi/2}]$ and $w$ does not depend on $\phi$, a uniformly attaining measurement exists if and only if the model is a great circle on the sphere and the equator defines the planes of the uniformly attaining measurements. In the same way, if $\phi$ is known and $\eta$ is the unknown parameter of interest, then uniform attainability can be achieved if and only if the model is of constant longitude, through the north and south poles. Under certain regularity conditions, models that admit uniformly attaining measurements belong to the class of quantum exponential transformation models [see Amari and Nagaoka (2000)].

## 5. Proofs.

PROOF OF LEMMA 1 (The SLD of a mixed state). Replacing the expressions for $\rho(\theta)$ and $\rho_{/\theta}$ in (3) gives (all arguments omitted)

$$
\begin{aligned}
(10) \quad & w_{/\theta}(\rho_1 - \rho_2) + w\rho_{1/\theta} + (1-w)\rho_{2/\theta} \\
& = \tfrac{1}{2}\{[w\rho_1 + (1-w)\rho_2]\rho_{//\theta} + \rho_{//\theta}[w\rho_1 + (1-w)\rho_2]\}.
\end{aligned}
$$

Pre- and postmultiplying both the members of (10) by $\langle\psi_h(\theta)|$ and $|\psi_k(\theta)\rangle$, $h, k = 1, 2$, respectively, gives, by (viii), $[\rho^\Psi_{//\theta}]_{11} = w_{/\theta}w^{-1}$, $[\rho^\Psi_{//\theta}]_{22} = -w_{/\theta}(1-w)^{-1}$ and $[\rho^\Psi_{//\theta}]_{12} = (2w-1)I_1^{1/2} = [\rho^\Psi_{//\theta}]_{21}$, where $[\rho^\Psi_{//\theta}]_{hk} = \langle\psi_h|\rho_{k/\theta}|\psi_h(\theta)\rangle$ indicates the $hk$th element of the matrix $\rho_{//\theta}$ with respect to the ordered basis $\Psi = \{|\psi_1\rangle, |\psi_2\rangle\}$. Hence, in the original coordinate system, $\rho_{//\theta} = (w_{/\theta}w^{-1})\rho_1 + (2w-1)I_1^{1/2}(\rho_{12} + \rho_{21}) - w_{/\theta}(1-w)^{-1}\rho_2$ where $\rho_{12} = |\psi_1(\theta)\rangle\langle\psi_2(\theta)|$ and $\rho_{21} = |\psi_2(\theta)\rangle\langle\psi_1(\theta)|$. To complete the proof, note that

$$
\begin{aligned}
I_1^{1/2}(\rho_{12} + \rho_{21}) &= I_1^{1/2}(|\psi_1(\theta)\rangle\langle\psi_1(\theta)\rho_{1//\theta}|I_1^{-1/2} + I_1^{-1/2}|\rho_{1//\theta}\psi_1(\theta)\rangle\langle\psi_1(\theta)|) \\
&= 2\rho_{1/\theta} = \rho_{1//\theta}.
\end{aligned}
$$

Conversely, if (6) holds, then (10) is an identity, by the uniqueness of $\rho_{//\theta}$. □

PROOF OF LEMMA 2 (Quantum information in a mixed state). Replacing, in (2), $\rho(\theta)$ by (5) and $\rho^2_{//\theta}$ by (6), $\rho_1$ and $\rho_2$ being pure and orthogonal, by linearity of the trace operator and by (ii), (iii) and (vi)–(ix),

$$
\begin{aligned}
I(\theta) &= \frac{(w_{/\theta})^2}{w} + w(2w-1)^2 I_1 + (1-w)(2w-1)^2 I_1 + \frac{(w_{/\theta})^2}{1-w} \\
&= w_{/\theta}^2\left(\frac{1}{w} + \frac{1}{1-w}\right) + (2w-1)^2 I_1.
\end{aligned}
$$
□

PROOF OF THEOREM 1 (Attainability in a mixed state). If equality holds, then replacing in (4) $\rho(\theta)^{1/2}$ and $\rho_{//\theta}$ by their expressions derived from (5) and (6) and postmultiplying both the members by $|\psi_1(\theta)\rangle$, we get

$$
(11) \quad m(x)^{1/2}\left\{\left[k(x,\theta)^{1/2} - \frac{w_{/\theta}}{w}\right]|\psi_1(\theta)\rangle - (2w-1)I_1^{1/2}|\psi_2(\theta)\rangle\right\} = 0.
$$

Since $|\psi_1(\theta)\rangle$ and $|\psi_2(\theta)\rangle$ are orthogonal, for $w \neq \tfrac{1}{2}$, the above system makes sense if and only if the matrix $m^{1/2}(x)$ [and consequently $m(x)$] is singular. A $2 \times 2$ Hermitian, singular and nonnegative matrix is necessarily of the form

$$
(12) \quad m(x)^{1/2} = c(x)|\gamma(x)\rangle\langle\gamma(x)|,
$$

where $|\gamma(x)\rangle = \langle\gamma^*(x)|\,|\gamma^*(x)\rangle^{-1/2}|\gamma^*(x)\rangle$, $\langle\gamma^*(x)| = [1\ \overline{\alpha_x}]$, $c(x) = a^2(x)\langle\gamma^*(x)|\,|\gamma^*(x)\rangle$ with $a(x) \in \mathbb{R}$ and $\alpha_x \in \mathbb{C}$. Replacing expression (12) for $m(x)^{1/2}$ in the equality condition (11) and premultiplying both the members by $c^{-1}(x)\langle\gamma(x)|$, we obtain

$$(13) \qquad \langle\gamma(x)|\,|\psi_1(\theta)\rangle = r(x,\theta)\langle\gamma(x)|\,|\psi_2(\theta)\rangle$$

with $r(x,\theta) = [(2w-1)I_1^{1/2}][k(x,\theta)^{1/2} - \frac{w_{/\theta}}{w}]^{-1} \in \mathbb{R}$.

Postmultiplying both the members of (4) by $|\psi_2(\theta)\rangle$ instead of $|\psi_1(\theta)\rangle$, one gets

$$(14) \qquad \langle\gamma(x)|\,|\psi_1(\theta)\rangle = r'(x,\theta)\langle\gamma(x)|\,|\psi_2(\theta)\rangle,$$

where now $r'(x,\theta) = [k(x,\theta)^{1/2} + \frac{w_{/\theta}}{1-w}][(2w-1)I_2^{1/2}]^{-1}$.

Combining (13) and (14) gives $r(x,\theta) = r'(x,\theta)$, that is,

$$k(x,\theta) + k(x,\theta)^{1/2}\left[\frac{w_{/\theta}(2w-1)}{w(1-w)}\right] - I(\theta) = 0.$$

Solving the equation for $k(x,\theta)^{1/2}$, one obtains

$$k(\theta) = \frac{(w_{/\theta})^2(2w-1)^2}{2w^2(1-w)^2} + I(\theta) \mp \frac{w_{/\theta}(2w-1)}{w(1-w)}\left[\frac{1}{4}\frac{(w_{/\theta})^2(2w-1)^2}{w^2(1-w)^2} + I(\theta)\right]^{1/2}$$

which, replaced in $r(x,\theta)$ or $r'(x,\theta)$, gives the proportionality constant between $\langle\gamma(x)|\,|\psi_1(\theta)\rangle$ and $\langle\gamma(x)|\,|\psi_2(\theta)\rangle$.

On the other hand, if $m(x) = c(x)^2|\gamma(x)\rangle\langle\gamma(x)|$, $c(x) \in \mathbb{R}$ and (13) holds, then premultiplying both its members by $c(x)|\gamma(x)\rangle$, replacing $|\psi_2(\theta)\rangle$ by its expression as a function of $\rho_{1//\theta}$, writing the latter as derived by (6) and multiplying by $w^{1/2}(\theta)$ gives

$$(15) \qquad k(\theta)^{1/2}m(x)^{1/2}\rho(\theta)^{1/2}|\psi_1(\theta)\rangle = m(x)^{1/2}\rho_{//\theta}\rho(\theta)^{1/2}|\psi_1(\theta)\rangle,$$

since $k(x,\theta)^{1/2} = r(x,\theta)^{-1}I_1^{1/2}(\theta)(2w-1) + \frac{w_{/\theta}}{w}$.

Replacing $|\psi_1(\theta)\rangle$, instead of $|\psi_2(\theta)\rangle$, with its expression as a function of $\rho_{2//\theta}$ gives

$$(16) \qquad m(x)^{1/2}\rho_{//\theta}\rho(\theta)^{1/2}|\psi_2(\theta)\rangle = k^{1/2}(\theta)m(x)^{1/2}\rho(\theta)^{1/2}|\psi_2(\theta)\rangle,$$

where $k^{1/2}(x,\theta) = r'(x,\theta)(2w-1)I_2^{1/2}(\theta) - \frac{w_{/\theta}}{1-w}$.

To complete the proof, it has to be shown that the vector equalities (15) and (16) imply the (matrix) equality condition (4). It follows from (15) and (16) that both $|\psi_1(\theta)\rangle$ and $|\psi_2(\theta)\rangle$ belong to the null space of the matrix $k(x,\theta)^{1/2}m(x)^{1/2} \times \rho(\theta)^{1/2} - m(x)^{1/2}\rho_{//\theta}\rho(\theta)^{1/2}$. However, $|\psi_1(\theta)\rangle$ and $|\psi_2(\theta)\rangle$ constitute an orthogonal basis of $\mathbb{C}^2$, and therefore such a matrix must necessarily be the null matrix; that is to say, $k(x,\theta)^{1/2}m(x)^{1/2}\rho(\theta)^{1/2} = m(x)^{1/2}\rho_{//\theta}\rho(\theta)^{1/2}$. $\square$

**Acknowledgments.** I thank Ole Barndorff-Nielsen, Richard Gill and Peter Jupp kindly for useful discussions. I give special thanks also to two anonymous referees for their valuable comments.

## REFERENCES

Amari, S.-I. and Nagaoka, H. (2000). *Methods of Information Geometry*. Oxford Univ. Press. MR1800071

Barndorff-Nielsen, O. E. and Gill, R. D. (2000). Fisher information in quantum statistics. *J. Phys. A* **33** 4481–4490. MR1768745

Braunstein, S. L. and Caves, C. M. (1994). Statistical distance and the geometry of quantum states. *Phys. Rev. Lett.* **72** 3439–3443. MR1274631

Fujiwara, A. and Nagaoka, H. (1995). Quantum Fisher metric and estimation for pure state models. *Phys. Lett. A* **201** 119–124. MR1329961

Helstrom, C. W. (1967). Minimum mean-squared error of estimates in quantum statistics. *Phys. Lett. A* **25** 101–102.

Helstrom, C. W. (1976). *Quantum Detection and Estimation Theory*. Academic Press, New York.

Holevo, A. S. (1982). *Probabilistic and Statistical Aspects of Quantum Theory*. North-Holland, Amsterdam. MR681693

Luati, A. (2003). A real formula for transition probability. *Statistica* **63**. To appear.

Department of Statistics
University of Bologna
Via Belle Arti 41
40126 Bologna
Italy
e-mail: luati@stat.unibo.it